\documentclass[a4paper,leqno,twoside]{article}
\usepackage{latexsym,amssymb,amsmath, amsthm}
\usepackage[affil-it]{authblk}
\newtheorem{theo}{Theorem}[section]
\newtheorem{coro}[theo]{Corollary}
\newtheorem{lemm}[theo]{Lemma}
\newtheorem{prop}[theo]{Proposition}

\newtheorem{rema}[theo]{Remark}

\author{I.D. Chipchakov}
\affil{Institute of Mathematics and Informatics\\ Bulgarian Academy
of Sciences\\ Acad. G. Bonchev Str., bl. 8\\ 1113, Sofia, Bulgaria;
email: chipchak@math.bas.bg}

\begin{document}
\title{On the behaviour of Brauer $p$-dimensions under
finitely-generated field extensions\footnote{Throughout this paper,
we write for brevity "FG-extension(s)" instead of "finitely-generated [field] extension(s)".}}

\date{\today}
\maketitle

\begin{abstract}
The present paper shows that if $q \in \mathbb P$ or $q = 0$, where
$\mathbb P$ is the set of prime numbers, then there exist
characteristic $q$ fields $E _{q,k}\colon \ k \in \mathbb N$, of
Brauer dimension Brd$(E _{q,k}) = k$ and infinite absolute Brauer
$p$-dimensions abrd$_{p}(E _{q,k})$, for all $p \in \mathbb P$ not 
dividing $q ^{2} - q$. This ensures that Brd$_{p}(F _{q,k}) =
\infty $, $p \dagger q ^{2} - q$, for every finitely-generated 
transcendental extension $F _{q,k}/E _{q,k}$. We also prove that
each sequence $a _{p}, b _{p}$, $p \in \mathbb P$, satisfying the
conditions $a _{2} = b _{2}$ and $0 \le b _{p} \le a _{p} \le
\infty $, equals the sequence abrd$_{p}(E), {\rm Brd}_{p}(E)$, $p
\in \mathbb P$, for a field $E$ of characteristic zero.
\end{abstract}

(secondary).

\par
\section{\bf Introduction}
\par
\medskip
Let $E$ be a field, $s(E)$ the class of finite-dimensional
associative central simple $E$-algebras, $d(E)$ the subclass of
division algebras $D \in s(E)$, and for each $A \in s(E)$, let $[A]$
be the equivalence class of $A$ in the Brauer group Br$(E)$. It is
known that Br$(E)$ is an abelian torsion group (cf. \cite{P}, Sect.
14.4), whence it decomposes into the direct sum of its
$p$-components Br$(E) _{p}$, where $p$ runs across the set $\mathbb
P$ of prime numbers. By Wedderburn's structure theorem (see, e.g.,
\cite{P}, Sect. 3.5), each $A \in s(E)$ is isomorphic to the full
matrix ring $M _{n}(D _{A})$ of order $n$ over some $D _{A} \in
d(E)$, uniquely determined by $A$, up-to an $E$-isomorphism. This
implies the dimension $[A\colon E]$ is a square of a positive
integer deg$(A)$, the degree of $A$. The main numerical invariants of
$A$ are deg$(A)$, the Schur index ind$(A) = {\rm deg}(D _{A})$, and
the exponent exp$(A)$, i.e. the order of $[A]$ in Br$(E)$. The
following statements describe basic divisibility relations between 
ind$(A)$ and exp$(A)$, and give an idea of their behaviour under the 
scalar extension map Br$(E) \to {\rm Br}(R)$, in case $R/E$ is a
field extension of finite degree $[R\colon E]$ (see, e.g., \cite{P}, 
Sects. 13.4, 14.4 and 15.2):
\par
\medskip\noindent
(1.1) (a) $({\rm ind}(A), {\rm exp}(A))$ is a Brauer pair, i.e.
exp$(A)$ divides ind$(A)$ and is divisible by every $p \in \mathbb
P$ dividing ind$(A)$.
\par
(b) ind$(A \otimes _{E} B) = {\rm ind}(A){\rm ind}(B)$, if $B \in
s(E)$ and g.c.d.$\{{\rm ind}(A), {\rm ind}(B)\}$ $= 1$; in this case,
if $A, B \in d(E)$, then the tensor product $A \otimes _{E} B$ lies
in $d(E)$.
\par
(c) ind$(A)$, ind$(A \otimes _{E} R)$, exp$(A)$ and exp$(A \otimes
_{E} R)$ divide ind$(A \otimes _{E} R)[R\colon E]$, ind$(A)$, exp$(A
\otimes _{E} R)[R\colon E]$ and exp$(A)$, respectively.
\par
\medskip
Statements (1.1) (a), (b) imply Brauer's Primary Tensor Product
Decomposition Theorem, for any $\Delta \in d(E)$ (cf. \cite{P},
Sect. 14.4). Also, (1.1) (a) fully describes general restrictions on
index-exponent pairs, in the following sense:
\par
\medskip
(1.2) Given a Brauer pair $(m ^{\prime }, m) \in \mathbb N ^{2}$,
there is a field $F$ with $({\rm ind}(D), {\rm exp}(D))$ $= (m
^{\prime }, m)$, for some $D \in d(F)$ (Brauer, see \cite{P}, Sect.
19.6). One may take as $F$ any rational (i.e. purely transcendental)
extension in infinitely many variables over any fixed field $F _{0}$.
\par
\medskip
The Brauer $p$-dimensions Brd$_{p}(E)$, $p \in \mathbb P$, of a
field $E$ contain essential information about pairs ind$(D), {\rm
exp}(D)$, $D \in d(E)$. We say that Brd$_{p}(E)$ is finite and equal
to $n$, for a fixed $p \in \mathbb P$, if $n$ is the least integer
$\ge 0$, for which ind$(D _{p}) \le {\rm exp}(D _{p}) ^{n}$ whenever
$D _{p} \in d(E)$ and $[D _{p}] \in {\rm Br}(E) _{p}$. If no such
$n$ exists, we set Brd$_{p}(E) = \infty $. The absolute Brauer
$p$-dimension of $E$ is defined as the supremum abrd$_{p}(E) = {\rm
sup}\{{\rm Brd}_{p}(R)\colon \ R \in {\rm Fe}(E)\}$, where Fe$(E)$
is the set of finite extensions of $E$ in a separable closure $E
_{\rm sep}$. We have abrd$_{p}(E) = 0$, for some $p \in \mathbb P$,
$p \neq {\rm char}(E)$, if and only if the absolute Galois group
$\mathcal{G}_{E} = \mathcal{G}(E _{\rm sep}/E)$ is of cohomological
$p$-dimension cd$_{p}(\mathcal{G}_{E}) \le 1$ (cf. \cite{Se}, Ch.
II, 3.1). When $E$ is virtually perfect, i.e. char$(E) = 0$ or
char$(E) = q > 0$ and $E$ is a finite extension of its subfield $E
^{q} = \{e ^{q}\colon e \in E\}$, the following holds:
\par
\medskip
(1.3) Brd$_{p}(E ^{\prime }) \le {\rm abrd}_{p}(E)$, for all $p \in
\mathbb P$ and finite extensions $E ^{\prime }/E$.
\par
\medskip\noindent
The assertion is obvious, if char$(E) = 0$. If char$(E) = q
> 0$, then $[E ^{\prime }\colon E ^{\prime q}] = [E\colon E ^{q}]$,
for every finite extension $E ^{\prime }/E$ (cf. \cite{L}, Ch. VII,
Sect. 7). Therefore, (1.3) can be deduced from (1.1) (c) and
Albert's theory of $q$-algebras \cite{A1}, Ch. VII, Theorem~28 (see
also Lemma \ref{lemm4.1}).
\par
\medskip
It is known that ${\rm Brd}_{p}(E) = {\rm abrd}_{p}(E) = 1$, for all
$p \in \mathbb P$, if $E$ is a global or local field (cf. \cite{Re},
(31.4) and (32.19)), or the function field of an algebraic surface
defined over an algebraically closed field $E _{0}$ \cite{Jong},
\cite{Lieb}. As shown in \cite{Matz}, when $E$ is the function field
of an $n$-dimensional algebraic variety over the field $E _{0}$, we
have abrd$_{p}(E) < p ^{n-1}$, $p \in \mathbb P$. The suprema
Brd$(E) = {\rm sup}\{{\rm Brd}_{p}(E)\colon \ p \in \mathbb P\}$ and
abrd$(E) = {\rm sup}\{{\rm Brd}(R)\colon \ R \in {\rm Fe}(E)\}$ are
called a Brauer dimension and an absolute Brauer dimension of $E$,
respectively. In view of (1.1), the definition of Brd$(E)$ is the
same as in \cite{ABGV}, Sect. 4. It has recently been proved
\cite{HHKr}, \cite{PaSur} (see also Lemmas \ref{lemm4.3} and
\ref{lemm4.4}), that abrd$(K _{m}) < \infty $, if $(K _{m}, v _{m})$
is an $m$-dimensional local field, in the sense of \cite{FV}, with a
quasifinite ($m$-th) residue field.
\par
\medskip
The present research considers the sequence Brd$_{p}(F)$, $p \in
\mathbb P$, for a transcendental FG-extension $F$ of a field $E$,
and its dependence upon abrd$_{p}(E)$, $p \in \mathbb P$. It is
motivated mainly by an open problem posed in Section 4 of the
survey \cite{ABGV}.

\medskip
\section{\bf The main results}
\par
\medskip
Fields $E$ with abrd$_{p}(E) < \infty $, for all $p \in \mathbb P$,
are singled out by Galois cohomology (see Remark \ref{rema4.2}), and
in the virtually perfect case, by the validity of the Primary Tensor
Product Decomposition Theorem, for every locally finite-dimensional
associative central division $E$-algebra of at most countable
dimension (see (1.3) and \cite{Ch1}). The applicability of this
result to basic fields of algebraic number theory and algebraic
geometry raises interest in the open problem of whether
FG-extensions of a global, local or algebraically closed field are
of finite absolute Brauer dimensions. This draws our attention to
the following open question:
\par
\medskip\noindent
(2.1) Is the class of fields $E$ of finite absolute Brauer
$p$-dimensions, for a fixed $p \in \mathbb P$, $p \neq {\rm
char}(E)$, closed under the formation of FG-extensions?
\par
\medskip
The purpose of this paper is to answer the similar question of
whether Brd$(F) < \infty $, for every FG-extension $F$ of a field
$E$ with Brd$(E) < \infty $ (this is stated in \cite{ABGV} as
Problem~4.4). Our starting point is the following result of
\cite{Ch4}:
\par
\medskip
\begin{prop}
\label{prop2.1} Let $E$ be a field, $p \in \mathbb P$ and $F/E$ an
{\rm FG}-extension of transcendency degree {\rm trd}$(F/E) = \kappa
\ge 1$. Then:
\par
{\rm (a)} {\rm Brd}$_{p}(F) \ge {\rm abrd}_{p}(E) + \kappa - 1$, if
{\rm abrd}$_{p}(E) < \infty $ and $F/E$ is rational;
\par
{\rm (b)} If {\rm abrd}$_{p}(E) = \infty $, then {\rm Brd}$_{p}(F) =
\infty $, and for any $n, m \in \mathbb N$ with $n \ge m > 0$, there
is $D _{n,m} \in d(F)$, such that {\rm ind}$(D _{n,m}) = p ^{n}$ and
{\rm exp}$(D _{n,m}) = p ^{m}$;
\par
{\rm (c)} Brd$_{p}(F) = \infty $, provided that $p = {\rm char}(E)$
and $[E\colon E ^{p}] = \infty $.
\end{prop}
\par
\medskip
The main result of the present paper can be stated as follows:
\par
\medskip
\begin{theo}
\label{theo2.2} For each $q \in \mathbb P \cup \{0\}$ and $k \in
\mathbb N$, there exists a field $E _{q,k}$ with char$(E _{q,k}) =
q$, Brd$(E _{q,k}) = k$ and abrd$_{p}(E _{q,k}) = \infty $, for all
$p \in \mathbb P \setminus P _{q}$, where $P _{0} = \{2\}$ and $P
_{q} = \{p \in \mathbb P\colon \ p \mid q ^{2} - q)\}$, $q \in
\mathbb P$. Moreover, if $q > 0$, then $E _{q,k}$ can be chosen so
that $[E _{q,k}\colon E _{q,k} ^{q}] = \infty $.
\end{theo}
\par
\medskip
When $q = 0$, the assertion of Theorem \ref{theo2.2} is contained in
our next result, which clarifies with Proposition \ref{prop2.1} the
influence of invariants abrd$_{p}(E)$, $p \in \mathbb P$, on the
behaviour of Brd$_{p}(F)$, $p \in \mathbb P$, for any transcendental
FG-extension $F/E$:
\par
\medskip
\begin{theo}
\label{theo2.3} Let $b _{p}, a _{p}$, $p \in \mathbb P$, be a
sequence with terms in the set $\mathbb N _{\infty } = \mathbb N
\cup \{0, \infty \}$, such that $b _{2} = a _{2}$ and $b _{p} \le
a _{p} \le \infty $, $p \in \mathbb P$. Then there exists a Henselian
field $(K, v)$ with {\rm char}$(\widehat K) = 0$,
$\mathcal{G}_{\widehat K}$ pronilpotent of cohomological dimension
{\rm cd}$(\mathcal{G}_{\widehat K}) \le 1$, and $({\rm abrd}_{p}(K),
{\rm Brd}_{p}(K)) = (a _{p}, b _{p})$, $p \in \mathbb P$.
\end{theo}
\par
\medskip
Proposition \ref{prop2.1}, Theorem \ref{theo2.2} and statement
(1.1) (b) imply the following:
\par
\medskip\noindent
(2.2) There exist fields $E _{k}$, $k \in \mathbb N$, such that
char$(E _{k}) = 2$, Brd$(E _{k}) = k$ and all Brauer pairs $(m
^{\prime }, n ^{\prime }) \in \mathbb N ^{2}$ are index-exponent
pairs over any transcendental FG-extension of $E _{k}$.
\par
\medskip\noindent
It is not known whether (2.2) holds in any characteristic $q \neq
2$. This is closely related to the following open problem:
\par
\medskip\noindent
(2.3) Find whether there exists a field $E$ containing a primitive
$p$-th root of unity, for a given $p \in \mathbb P$, such that
Brd$_{p}(E) < {\rm abrd}_{p}(E) = \infty $.
\par
\medskip
Statement (1.1) (b), Proposition \ref{prop2.1} and Theorem
\ref{theo2.2} imply the validity of (2.2), for $q = 0$ and Brauer
pairs of odd positive integers. When $q > 2$, they show that if $[E _{q,k}\colon E _{q,k} ^{q}] = \infty $, then Brauer pairs $(m
^{\prime }, m) \in \mathbb N ^{2}$ relatively prime to $q - 1$ are
index-exponent pairs over every transcendental FG-extension of $E _{q,k}$. Thus Problem~4.4 of \cite{ABGV} is solved in the negative.
As a whole, our research shows that (2.1) can be a suitable
replacement in the list of \cite{ABGV} for this problem.
\par
\medskip
The proofs of our main results are based on results of valuation
theory like Morandi's theorem on tensor products of valued division
algebras \cite{Mo}, Theorem~1, and the classical Ostrowski theorem.
They rely on a standard method of realizing profinite groups as
Galois groups \cite{Wat}, and on a construction of Henselian fields
with prescribed properties of their value groups, residue fields and
finite extensions. We also use a characterization of fields $E$ with
abrd$_{p}(E) \le \mu $, for a given $\mu \in \mathbb N$ (which
generalizes Albert's theorem \cite{A1}, Ch. XI, Theorem~3), as well
as formulae for Brd$_{p}(K)$ and abrd$_{p}(K)$ concerning some
Henselian fields. This approach enables one to obtain the following:
\par
\medskip\noindent
(2.4) (a) There exists a field $E _{1}$ with abrd$(E _{1}) = \infty
$, abrd$_{p}(E _{1}) < \infty $, $p \in \mathbb P$, and Brd$(L _{1})
< \infty $, for every finite extension $L _{1}/E _{1}$; hence, by
\cite{Ch4}, Corollary~5.4, Brd$(F _{1}) = \infty $, for every transcendental FG-extension $F _{1}/E _{1}$;
\par
(b) For any integer $n \ge 2$, there is a Galois extension $L _{n}/E
_{n}$, such that $[L _{n}\colon E _{n}] = n$, Brd$_{p}(L _{n}) =
\infty $, for all $p \in \mathbb P$, $p \equiv 1 ({\rm mod} \ n)$,
and Brd$(M _{n}) < \infty $, provided that $M _{n}$ is an extension
of $E$ in $L _{n,{\rm sep}}$ not including $L _{n}$.
\par
\medskip
Our basic notation and terminology are standard, as used in
\cite{Ch2}. For any field $K$ with a Krull valuation $v$, unless
stated otherwise, we denote by $\widehat K$ and $v(K)$ the residue
field and the value group of $(K, v)$, respectively; $v(K)$ is
supposed to be an additively written totally ordered abelian group.
As usual, $\mathbb Z$ stands for the additive group of integers,
$\mathbb Z _{p}$ is the additive groups of $p$-adic integers, for
any $p \in \mathbb P$, and $[r]$ is the integral part of any real
number $r \ge 0$. We write $I(\Lambda ^{\prime }/\Lambda )$ for the
set of intermediate fields of a field extension $\Lambda ^{\prime
}/\Lambda $, and Br$(\Lambda ^{\prime }/\Lambda )$ for the relative
Brauer group of $\Lambda ^{\prime }/\Lambda $. By a $\Lambda
$-valuation of $\Lambda ^{\prime }$, we mean a Krull valuation $v$
with $v(\lambda ) = 0$, $\lambda \in \Lambda ^{\ast }$. Given a
field $E$ and $p \in \mathbb P$, $E(p)$ denotes the maximal
$p$-extension of $E$ in $E _{\rm sep}$, and $r _{p}(E)$ the rank
of the Galois group $\mathcal{G}(E(p)/E)$ as a pro-$p$-group ($r
_{p}(E) = 0$, if $E(p) = E$). Brauer groups are considered to be
additively written, Galois groups are viewed as profinite with
respect to the Krull topology, and by a homomorphism of profinite
groups, we mean a continuous one. We refer the reader to
\cite{Efr2}, \cite{JW}, \cite{L}, \cite{P} and \cite{Se}, for any
missing definitions concerning valuation theory, field extensions,
simple algebras, Brauer groups and Galois cohomology.
\par
\medskip
Here is an overview of the rest of the paper: Section 3 includes
preliminaries used in the sequel, and Galois-theoretic ingredients
of the proof of Theorem \ref{theo2.3}. Theorems \ref{theo2.2},
\ref{theo2.3} and statement (2.4) are proved in Section 4.

\medskip
\section{\bf Preliminaries on Henselian valuations and preparation for
the proof of Theorem \ref{theo2.3}}
\par
\medskip
The results of this Section are known and will often be used without
an explicit reference. Assume that $(K, v)$ is a Henselian field,
i.e. $v$ is a Krull valuation on $K$, which extends uniquely, up-to
an equivalence, to a valuation $v _{L}$ on each algebraic extension $L/K$. Put $v(L) = v _{L}(L)$ and denote by $\widehat L$ the
residue field of $(L, v _{L})$. It is known that $\widehat L/\widehat
K$ is an algebraic extension and $v(K)$ is a subgroup of $v(L)$. When $[L\colon K]$ is finite, Ostrowski's theorem states the following
(cf. \cite{Efr2}, Theorem~17.2.1):
\par
\medskip\noindent
(3.1) $[\widehat L\colon \widehat K]e(L/K)$ divides $[L\colon K]$ and
$[L\colon K][\widehat L\colon \widehat K] ^{-1}e(L/K) ^{-1}$ is not
divisible by any $p \in \mathbb P$ different from char$(\widehat K)$,
$e(L/K)$ being the index of $v(K)$ in $v(L)$;  in particular, if
char$(\widehat K) \dagger [L\colon K]$, then $[L\colon K] = [\widehat
L\colon \widehat K]e(L/K)$.
\par
\medskip
Statement (3.1) and the Henselity of $v$ imply the following:
\par
\medskip\noindent
(3.2) The quotient groups $v(K)/pv(K)$ and $v(L)/pv(L)$ are
isomorphic, if $p \in \mathbb P$ and $L/K$ is a finite extension.
When char$(\widehat K) \dagger [L\colon K]$, the natural embedding
of $K$ into $L$ induces canonically an isomorphism $v(K)/pv(K) \cong
v(L)/pv(L)$.
\par
\medskip
A finite extension $R/K$ is said to be defectless, if $[R\colon K] =
[\widehat R\colon \widehat K]e(R/K)$. It is called inertial, if
$[R\colon K] = [\widehat R\colon \widehat K]$ and $\widehat
R/\widehat K$ is separable. We say that $R/K$ is totally ramified,
if $[R\colon K] = e(R/K)$. The Henselity of $v$ ensures that the
compositum $K _{\rm ur}$ of inertial extensions of $K$ in $K _{\rm
sep}$ has the following properties:
\par
\medskip\noindent
(3.3) (a) $v(K _{\rm ur}) = v(K)$ and finite extensions of $K$ in $K
_{\rm ur}$ are inertial;
\par
(b) $K _{\rm ur}/K$ is a Galois extension, $\widehat K _{\rm ur}
\cong \widehat K _{\rm sep}$ over $\widehat K$, $\mathcal{G}(K _{\rm
ur}/K)$ $\cong \mathcal{G}_{\widehat K}$, and the natural mapping of
$I(K _{\rm ur}/K)$ into $I(\widehat K _{\rm sep}/\widehat K)$ is
bijective.
\par
\medskip\noindent
When $(K, v)$ is Henselian, each $\Delta \in d(K)$ has a unique,
up-to an equivalence, valuation $v _{\Delta }$ extending $v$ so that
the value group $v(\Delta )$ of $(\Delta , v _{\Delta })$ is totally
ordered and abelian (cf. \cite{Sch}, Ch. 2, Sect. 7). It is known
that $v(K)$ is a subgroup of $v(\Delta )$ of index $e(\Delta /K) \le
[\Delta \colon K]$, and the residue division ring $\widehat \Delta $
of $(\Delta , v _{\Delta })$ is a $\widehat K$-algebra. By the
Ostrowski-Draxl theorem \cite{Dr1}, $[\Delta \colon K]$ is divisible
by $e(\Delta /K)[\widehat {\Delta }\colon \widehat K]$, and in case
char$(\widehat K) \dagger [\Delta \colon K]$, $[\Delta \colon K] =
e(\Delta /K)[\widehat \Delta \colon \widehat K]$. An algebra $D \in
d(K)$ is called inertial, if $[D\colon K] = [\widehat D\colon
\widehat K]$ and $\widehat D \in d(\widehat K)$. Inertial
$K$-algebras and those in $d(\widehat K)$ are related as follows
(see \cite{JW}, Theorem~2.8):
\par
\medskip\noindent
(3.5) (a) Each $\widetilde D \in d(\widehat K)$ has an inertial lift
over $K$, i.e. $\widetilde D = \widehat D$, for some $D \in d(K)$
inertial over $K$, and uniquely determined by $\widetilde D$, up-to
a $K$-isomorphism.
\par
(b) The set IBr$(K) = \{[I] \in {\rm Br}(K)\colon \ I \in d(K)$ is
inertial$\}$ is a subgroup of Br$(K)$; the canonical mapping IBr$(K)
\to {\rm Br}(\widehat K)$ is an isomorphism.
\par
\medskip
The study of Brd$_{p}(K)$, for a given $p \in \mathbb P$, relies on
constructive methods based on the following statements:
\par
\medskip
(3.6) (a) If $U _{1}, \dots , U _{n}$ are cyclic $p$-extensions of
$K$ in $K _{\rm ur}$, and $\pi _{1}, \dots , \pi _{n}$ are elements
of $K ^{\ast }$, such that $[U _{1} \dots U _{n}\colon K] = \prod
_{j=1} ^{n} [U _{j}\colon K]$ and the cosets $\bar \pi _{j} = v(\pi
_{j}) + pv(K)$, $j = 1, \dots , n$, are linearly independent over
$\mathbb F _{p}$, then $d(K)$ contains the $K$-algebra $B = \otimes
_{j=1} ^{n} B _{j}$, where $\otimes = \otimes _{K}$, and for each
index $j$, $B _{j}$ is the cyclic $K$-algebra $(U _{j}/K, \tau _{j},
\pi _{j})$, $\tau _{j}$ being a generator of $\mathcal{G}(U _{j}/K)$.
\par
(b) If $p \neq {\rm char}(\widehat K)$, $K$ contains a primitive
$p$-th root of unity $\varepsilon $, and $\pi _{1} ^{\prime }, \dots
, \pi _{2m} ^{\prime }$ are elements of $K ^{\ast }$, such that
$\bar \pi _{i} ^{\prime } = v(\pi _{i} ^{\prime }) + pv(K)$, $i = 1,
\dots , 2m$, are linearly independent over $\mathbb F _{p}$, then the
$K$-algebra $T = \otimes _{u=1} ^{2m} T _{u}$ lies in $d(K)$, where
$\otimes = \otimes _{K}$ and $T _{u}$ is the symbol $K$-algebra $A
_{\varepsilon }(\pi _{2u-1} ^{\prime }, \pi _{2u} ^{\prime }; K)$,
for every index $u$.
\par
(c) Under the hypotheses of (a) and (b), if the system $\bar \pi
_{j}, \bar \pi _{i} ^{\prime }$, $j = 1, \dots , n$; $i = 1, \dots
2m$, is linearly independent over $\mathbb F _{p}$, then $B \otimes
_{K} T \in d(K)$.
\par
\medskip\noindent
Statement (3.6) (c) follows at once from \cite{Mo}, Theorem~1, and
(3.6) (a) is a special case of \cite{JW}, Example~4.3. Also, it is
clear from Kummer theory and the conditions of (3.6) (b) that $T
_{1}, \dots , T _{m}$ are cyclic $K$-algebras. Using (3.1) and the
Henselity of $v$, one obtains further that $v(\pi _{2u} ^{\prime })
\neq v(\lambda _{u})$, for any element $\lambda _{u}$ of the norm
group $N(K(\sqrt[p]{\pi _{2u-1} ^{\prime }})/K)$. Therefore, $\pi
_{2u} ^{\prime } \notin N(K(\sqrt[p]{\pi _{2u-1} ^{\prime }})/K)$,
so it follows from well-known general properties of cyclic
$K$-algebras (cf. \cite{P}, Sect. 15.1, Proposition~b) that $T _{u}
\in d(K)$, $u = 1, \dots , m$. It is now easily deduced from
\cite{Mo}, Theorem~1, that $T \in d(K)$, as claimed.

\medskip
Statements (3.6) and the following lemma play a crucial role in the
proof of Theorem \ref{theo2.2}.

\medskip
\begin{lemm}
\label{lemm3.1} Let $K _{0}$ be a perfect field with {\rm char}$(K
_{0}) = q \ge 0$, and let $n(p)\colon \ p \in \mathbb P$, be a
sequence with terms in $\mathbb N _{\infty }$. Then there exists a
Henselian field $(K, v)$, such that {\rm char}$(K) = q$, $\widehat K
= K _{0}$, and for each $p \in \mathbb P$, the group $v(K)/pv(K)$
has dimension $n(p)$ as a vector space over the field $\mathbb F
_{p}$ with $p$ elements. When $q > 0$ and $n(q) < \infty $, $K$ can
be chosen so that $[K\colon K ^{q}] = q ^{n(q)}$ and $r _{q}(K) =
\infty $, and in case $n(q) > 0$, $v(K)$ possesses an isolated
subgroup $H$ satisfying the following:
\par
{\rm (a)} $H/pH \cong v(K)/pv(K)$ and $v(K)/H = p(v(K)/H)$, $p \in
\mathbb P \setminus \{q\}$; $H = qH$ and $(v(K)/H)/q(v(K)/H) \cong
v(K)/qv(K)$;
\par
{\rm (b)} The valuation $v _{H}$ of $K$ with $v _{H}(K) = v(K)/H$,
defined by the composition $\eta _{H} \circ v\colon \ K ^{\ast } \to
v(K)/H$, where $\eta _{H}$ is the natural homomorphism of $v(K)$
upon $v(K)/H$, has a perfect residue field $K _{H} \in I(K/K _{0})$
with $r _{q}(K _{H}) = \infty $.
\end{lemm}

\medskip
\begin{proof}
Let $K _{\infty }$ be an extension of $K _{0}$ obtained as the union
$K _{\infty } = \cup _{n \in \mathbb N} K _{n}$ of iterated formal
(Laurent) power series fields, defined inductively by the rule $K
_{n} = K _{n-1}((X _{n}))$, $n \in \mathbb N$. Denote by $\omega
_{n}$ the standard $K _{0}$-valuation of $K _{n}$ with $\omega
_{n}(K _{n}) = \mathbb Z ^{n}$, for each $n \in \mathbb N$ ($\mathbb
Z ^{n}$ is viewed as an ordered group with respect to the inverse
lexicographic ordering $\le _{n}$). Let $\overline K _{\infty }$ be
an algebraic closure of $K _{\infty }$, $\omega $ the natural
valuation of $K _{\infty }$ extending $\omega _{n}$, for every $n$,
and in case $q > 0$, let $K _{\infty } ^{\prime }$ be the perfect
closure of $K _{\infty }$ in $\overline K _{\infty }$, and $\omega
^{\prime }$ a valuation of $K _{\infty } ^{\prime }$ extending
$\omega $. Clearly, $K _{0}$ is the residue field of $(K _{\infty }, 
\omega )$ and $\omega (K _{\infty })$ equals the union $\mathbb Z 
^{\infty } = \cup _{n \in \mathbb N} Z ^{n}$, considered with its
unique ordering inducing $\le _{n}$ on $\mathbb Z ^{n}$, for all $n$.
It is known (cf. \cite{Efr2}, Sects. 4.2 and 18.4) that the
valuations $\omega _{n}$, $n \in \mathbb N$, are Henselian, which
implies $\omega $ is of the same kind. Note further that if $q > 0$,
then the set $\rho (K _{\infty }) = \{u ^{q} - u\colon u \in K
_{\infty }\}$ is a vector subspace of $K _{\infty }$ over its prime
subfield $\mathbb F$, and $\omega (u ^{q} - u) \in q\omega (K
_{\infty })$ whenever $\omega (u) < 0$. This implies that, for any
$\pi \in K _{\infty }$ with $\omega (\pi ) < 0$ and $\omega (\pi )
\notin q\omega (K _{\infty })$, the cosets $\pi ^{1+qm} + \rho (K
_{\infty })$, $m \in \mathbb N$, are linearly independent over
$\mathbb F$, so it follows from the Artin-Schreier theorem (cf.
\cite{L}, Ch. VIII, Sect. 6) that $r _{q}(K _{\infty }) = \infty $.
These observations show that if $n(p) = \infty $, for all $p \in
\mathbb P$, then it suffices for the proof of Lemma \ref{lemm3.1} to
put $(K, v) = (K _{\infty }, \omega )$. When $q > 0 = n(q)$ and
$n(p) = \infty $, $p \neq q$, one may take as $(K, v)$ the valued
field $(K _{\infty } ^{\prime }, \omega ^{\prime })$; in this case,
we put $\Theta _{0} = K _{\infty } ^{\prime }$.
\par
Assume now that the set $P = \{p \in \mathbb P \setminus
\{q\}\colon \ n(p) < \infty \}$ is nonempty, $\bar \omega (K
_{\infty })$ is a divisible hull of $\omega (K _{\infty })$, and for
each $R \in I(\overline K _{\infty }/K _{\infty })$, $\omega _{R}$
is the valuation of $R$ extending $\omega $ so that $\omega (R) =
\omega _{R}(R)$ be an ordered subgroup of $\bar \omega (K _{\infty
})$. For any $p \in P$ and each index $n > n(p)$, let $\Sigma _{p,n}
= \{Y _{p,n,m}\colon \ m \in \mathbb N\}$ be a subset of $\overline
K _{\infty }$, such that $Y _{p,n,1} ^{p} = X _{n}$ and $Y _{p,n,m}
^{p} = Y _{p,n,(m-1)}$, $m \ge 2$. Put $\Sigma = \cup _{p \in P}
\Sigma _{p}$, where $\Sigma _{p} = \cup _{n=n(p)+1} ^{\infty }
\Sigma _{p,n}$, for each $p \in P$, and denote by $\widetilde K$ the
extension of $K _{\infty }$ generated by $\Sigma $. It is easily
verified that finite extensions of $K _{\infty }$ in $\widetilde K$
are totally ramified, and for each $p \in \mathbb P \setminus
\{q\}$, $n(p)$ equals the dimension of $\omega (\widetilde
K)/p\omega (\widetilde K)$ as an $\mathbb F _{p}$-vector space. As
$r _{q}(K _{\infty }) = \infty $ in case $q > 0$, then $r
_{q}(\widetilde K) = \infty $ as well. These observations show that
$(\widetilde K, \omega _{\widetilde K})$ has the property required
by Lemma \ref{lemm3.1}, if $q = 0$ or $q > 0$ and $n(q) = \infty $.
Suppose that $q > 0$, $n(q) < \infty $, and denote by $\Theta _{0}$
the perfect closure of $\widetilde K$ in $\overline K _{\infty }$. As
$K _{0}$ is perfect, the basic theory of algebraic extensions (cf.
\cite{L}, Ch. VII, Proposition~12) implies that $\omega (\Theta
_{0}) = q\omega (\Theta _{0})$, and for each $p \in \mathbb P
\setminus \{q\}$, $\omega (\Theta _{0})/p\omega (\Theta _{0})$ has
dimension $n(p)$ over $\mathbb F _{p}$. Thus Lemma \ref{lemm3.1} is
proved in the case where $n(q) = 0$.
\par
It remains to consider the case of $0 < n(q) < \infty $. Let $H = 
\omega (\Theta _{0}) = \omega _{\Theta _{0}}(\Theta _{0})$, $n(q) =
n$, $\Theta _{n}$ be an iterated formal power series field in $n$
variables over $\Theta _{0}$, $\kappa $ the standard $\mathbb Z
^{n}$-valued $\Theta _{0}$-valuation of $\Theta _{n}$, and $w$ the
valuation of $\Theta _{n}$ extending $\omega _{\Theta _{0}}$ so that 
$H$ be an isolated subgroup of $w(\Theta _{n})$, $w(\Theta _{n})$ 
the direct sum $H \oplus \kappa (\Theta _{n})$, and $\kappa $ be 
induced canonically by $w$ and $H$ (cf. \cite{Efr2}, Sect. 4.2). 
Then \cite{Efr2}, Theorem~18.1.2, and \cite{W}, Theorem~32.15, imply 
$w$ inherits the Henselity of $\omega $ and $\kappa $. Applying 
(3.1), \cite{Efr2}, Theorem~18.4.1, and \cite{W}, Theorem~31.21, one
concludes that finite extensions of $\Theta _{n}$ are defectless
relative to $\kappa $, and $n$ equals the $\mathbb F _{p}$-dimension
of $\kappa (\Theta _{n})/p\kappa (\Theta _{n})$, for $p \in \mathbb
P$. Let now $K$ be a maximal extension of $\Theta _{n}$ in $\Theta
_{n,{\rm sep}}$ with respect to the property that finite extensions
of $\Theta _{n}$ in $K$ have degrees not divisible by $q$ and are
totally ramified over $\Theta _{n}$ relative to $\kappa $. Then
$[K\colon K ^{q}] = q ^{n}$, $\kappa (K) = p\kappa (K)$, $p \in
\mathbb P \setminus \{q\}$, and it follows from (3.2), \cite{Ch2},
(1.3), and the preceding observation that the natural embedding of
$\Theta _{n}$ into $K$ induces an isomorphism $\kappa (\Theta
_{n})/q\kappa (\Theta _{n}) \cong \kappa (K)/q\kappa (K)$. These
results and the obtained properties of $(\Theta _{0}, \omega
_{\Theta _{0}})$ indicate that $\kappa (K) \cong v(K)/H$ and $v(K)$ 
has the properties required by Lemma \ref{lemm3.1}, where $v = w 
_{K}$.
\end{proof}

\medskip
\begin{rema}
\label{rema3.2} Under the hypotheses of Lemma \ref{lemm3.1}, suppose
that $q > 0$ and $0 < n(q) = n < \infty $, $\Theta _{0}$, $\Theta 
_{n}$, $\kappa $, $w$ and $\omega $ are defined as in the proof of
the lemma, $\Theta _{j} = \Theta _{j-1}((Z _{j}))$, $j = 1, \dots ,
n$, and $\theta $ is the Henselian discrete $\Theta _{n-1}$-valuation
of $\Theta _{n}$. Denote by $\kappa '$, $w'$ and $\theta '$ the 
valuations of $K _{\rm sep}$ extending $\kappa $, $w$ and $\theta $, 
respectively, put $\Lambda _{0} = \Theta _{n-1}(Z _{n})$, and for any 
$\Lambda \in I(K _{\rm sep}/\Lambda _{0})$, let $w _{\Lambda }$, 
$\kappa _{\Lambda }$ and $\theta _{\Lambda }$ be the valuations of
$\Lambda $ induced by $w'$, $\kappa '$ and $\theta '$, respectively. 
Analyzing the proof of the latter part of Lemma \ref{lemm3.1}, one 
obtains the existence of a subset $\Gamma \subset K$, such that 
$\Theta _{n}(\Gamma ) = K$, the field $\Phi _{0} = \Lambda
_{0}(\Gamma )$ is separable over $\Lambda _{0}$, $\Phi _{0} ^{\prime
} \cap \Theta _{n} = \Lambda _{0}$, where $\Phi _{0} ^{\prime } \in
I(K _{\rm sep}/\Phi _{0})$ is the Galois closure of $\Phi _{0}$ over 
$\Lambda _{0}$, and for each $\Lambda \in I(\Theta _{n}/\Lambda
_{0})$, the field $\Phi = \Lambda (\Gamma )$ satisfies the condition 
$\kappa _{\Phi }(\Phi ) = \kappa (K)$, $\theta _{\Phi }(\Phi ) =
\theta (K)$ and $w _{\Phi }(\Phi ) = w(K)$. Also, $\Phi _{0} 
^{\prime }\Lambda $ is the root field of $\Phi $ in $K _{\rm sep}$ 
over $\Lambda $, $\mathcal{G}(\Phi _{0} ^{\prime }\Lambda /\Lambda ) 
\cong \mathcal{G}(\Phi _{0} ^{\prime }/\Lambda _{0})$ and finite 
extensions of $\Lambda $ in $\Phi $ are totally ramified of degrees 
not divisible by $q$. In addition, the residue fields of $w _{\Phi 
}$, $\kappa _{\Phi }$ and $\theta _{\Phi }$ are isomorphic to $K 
_{0}$, $\Theta _{0}$ and $\Theta _{n-1}$, respectively.
\end{rema}

\medskip
We conclude this Section with two lemmas which contain the main
Galois-theoretic ingredients of our proofs of (2.4) (a) and Theorem
\ref{theo2.3}. The former lemma makes it easy to prove Theorem
\ref{theo2.3} steering clear of (2.3).

\medskip
\begin{lemm}
\label{lemm3.3} There exists a field $E _{0}$ with {\rm char}$(E
_{0}) = 0$, such that $\mathcal{G}_{E _{0}}$ is isomorphic to the
additive group $\mathbb Z _{2}$ of $2$-adic integers, and for each
$p \in \mathbb P$, $[E _{0}(\varepsilon _{p})\colon E _{0}] = 2
^{y(p)}$, where $\varepsilon _{p}$ is a primitive $p$-th root of
unity in $E _{0,{\rm sep}}$, and $y(p)$ is the greatest integer for
which $2 ^{y(p)} \mid p - 1$.
\end{lemm}

\medskip
\begin{proof}
Let $\varepsilon _{p}$ be a primitive $p$-th root of unity in
$\mathbb Q _{\rm sep}$, and let $R _{0}$ be the extension of
$\mathbb Q$ in $\mathbb Q _{\rm sep}$ generated by the set $\Sigma =
\{\sqrt{\alpha _{p}}\colon p \in \mathbb P\}$, where $\alpha _{2} =
-2$ and $\alpha _{p} = (-1) ^{(p+1)/2}$, for each $p > 2$. Then it
follows from Kummer theory that $\sqrt{-1} \notin R _{0}$, i.e. the
set $\Sigma ^{\prime } = \{R \in I(\mathbb Q _{\rm sep}/R
_{0})\colon \sqrt{-1} \notin R\}$ is nonempty. Clearly, $\Sigma ^{\prime }$
satisfies the conditions of Zorn's lemma, whence it contains a
maximal element $E _{0}$ with respect to the partial ordering by
inclusion. In view of Galois theory, this ensures that Fe$(E _{0})$
consists of cyclic $2$-extensions. Observing also that $E _{0}$ is a
nonreal field (since $\sqrt{\alpha _{2}} \in E _{0}$), one obtains
from \cite{Wh}, Theorem~2, that $\mathcal{G}_{E _{0}} \cong \mathbb
Z _{2}$, as claimed. It remains to be seen that $[E _{0}(\varepsilon
_{p})\colon E _{0}] = 2 ^{y(p)}$, for an arbitrary fixed $p \in
\mathbb P \setminus \{2\}$. It is well-known that $\mathbb
Q(\varepsilon _{p})/\mathbb Q$ is a cyclic extension and $\mathbb
Q(\sqrt{\beta _{p}})$ is the unique quadratic extension of $\mathbb
Q$ in $\mathbb Q(\varepsilon _{p})$, where $\beta _{p} = (-1)
^{(p-1)/2}p$. It is therefore clear from Galois theory that the
equality $[E _{0}(\varepsilon _{p})\colon E _{0}] = 2 ^{y(p)}$ will
follow, if we show that $\sqrt{\beta _{p}} \notin E _{0}$. This,
however, is obvious, since $\beta _{p} = -\alpha _{p}$,
$\sqrt{\alpha _{p}} \in E _{0}$ and $\sqrt{-1} \notin E _{0}$, so
Lemma \ref{lemm3.3} is proved.
\end{proof}

\medskip
\begin{lemm}
\label{lemm3.4} Assume that $E _{0}$ is a field with {\rm
cd}$(\mathcal{G}_{E _{0}}) \le 1$, and $G$ is a profinite group with 
{\rm cd}$(G) \le 1$ and {\rm cd}$_{p}(G) = 0$ whenever $p \in 
\mathbb P$ and {\rm cd}$_{p}(\mathcal{G}_{E _{0}}) \neq 0$. Then
there exists a field extension $E/E _{0}$, such that $E _{0}$ is
algebraically closed in $E$ and $\mathcal{G}_{E}$ is isomorphic to
the topological group product $\mathcal{G}_{E _{0}} \times G$; in 
particular, all roots of unity in $E$ lie in $E _{0}$.
\end{lemm}

\medskip
\begin{proof}
It is known (cf. \cite{Wat}) that $E _{0}$ has extensions $R$ and $R
^{\prime }$, such that $R ^{\prime }/E _{0}$ is rational, $R \in I(R
^{\prime }/E _{0})$ and $R ^{\prime }/R$ is Galois with
$\mathcal{G}(R ^{\prime }/R) \cong G$. Identifying $E _{0,{\rm
sep}}$ with its $E _{0}$-isomorphic copy in $R ^{\prime }_{\rm
sep}$, and observing that $E _{0}$ is algebraically closed in $R
^{\prime }$, one obtains that $E _{0,{\rm sep}}R ^{\prime }/R$ is
Galois with $\mathcal{G}(E _{0,{\rm sep}}R ^{\prime }/R) \cong
\mathcal{G}_{E _{0}} \times G$. In view of the assumptions on
$\mathcal{G}_{E _{0}}$ and $G$, this yields cd$(\mathcal{G}(E
_{0,{\rm sep}}R ^{\prime }/R)) = 1$, which means that $\mathcal{G}(E
_{0,{\rm sep}}R ^{\prime }/R)$ is a projective profinite group (cf.
\cite{Se}, Ch. I, 5.9). Hence, by Galois theory, there is a field $E
\in I(R ^{\prime }_{\rm sep}/R)$, such that $E _{0,{\rm sep}}R
^{\prime }E = R ^{\prime }_{\rm sep}$ and $(E _{0,{\rm sep}}R
^{\prime }) \cap E = R$. This shows that $E _{0}$ is algebraically
closed in $E$ and $\mathcal{G}_{E} \cong \mathcal{G}(E _{0,{\rm
sep}}R ^{\prime }/R) \cong \mathcal{G}_{E _{0}} \times G$, which
proves Lemma \ref{lemm3.4}.
\end{proof}

\medskip
\section{\bf Proofs of Theorems 2.2 and 2.3}

\medskip
First we characterize the condition abrd$_{p}(E) \le \mu $,
for a field $E$ and a given $\mu \in \mathbb N$. When $E$ is
virtually perfect, by (1.3), this result in fact is equivalent to
\cite{PaSur}, Lemma~1.1, and in case $\mu = 1$, it restates
Theorem~3 of \cite{A1}, Ch. XI.

\medskip
\begin{lemm}
\label{lemm4.1} Let $E$ be a field, $p \in \mathbb P$ and $\mu \in
\mathbb N$. Then {\rm abrd}$_{p}(E) \le \mu $ if and only if, for
each $E ^{\prime } \in {\rm Fe}(E)$, {\rm ind}$(\Delta ) \le p ^{\mu
}$ whenever $\Delta \in d(E ^{\prime })$ and {\rm exp}$(\Delta ) =
p$.
\end{lemm}

\medskip
\begin{proof}
The left-to-right implication is obvious, so we prove only the
converse one. Fix a pair $E ^{\prime } \in {\rm Fe}(E)$, $\Delta
^{\prime } \in d(E ^{\prime })$ with exp$(\Delta ^{\prime }) = p
^{n}$, for some $n \in \mathbb N$. We show that ind$(\Delta ^{\prime
}) \mid p ^{n\mu }$. This is obvious, if $n = 1$, so we assume that
$n \ge 2$. Take $\Delta \in d(E ^{\prime })$ so that $[\Delta ] = p
^{n-1}[\Delta ^{\prime }]$, and let $Y$ be a maximal subfield of
$\Delta $. It is well-known that $[Y\colon E ^{\prime }] = {\rm
ind}(\Delta )$ and $Y$ can be chosen so as to be separable over $E
^{\prime }$ (see \cite{P}, Sect. 13.5). Therefore, our assumptions
show that $[Y\colon E ^{\prime }] \mid p ^{\mu }$. Note that, by the
choice of $\Delta $, $\Delta ^{\prime } \otimes _{E'} Y \in s(Y)$
and exp$(\Delta ^{\prime } \otimes _{E'} Y) = p ^{n-1}$. These
remarks and a standard inductive argument lead to the conclusion
that it suffices to prove the divisibility ind$(\Delta ^{\prime })
\mid p ^{n\mu }$, provided ind$(\Delta ^{\prime } \otimes _{E'} Y)
\mid p ^{(n-1)\mu }$. Fix $\Delta ^{\prime } _{Y} \in d(Y)$ so that
$[\Delta ^{\prime } _{Y}] = [\Delta ^{\prime } \otimes _{E'} Y]$,
and take a maximal subfield $Y ^{\prime }$ of $\Delta ^{\prime }
_{Y}$. Then $[Y ^{\prime }\colon E ^{\prime }] = {\rm ind}(\Delta
^{\prime } \otimes _{E'} Y).[Y\colon E ^{\prime }]$, which implies
$[Y ^{\prime }\colon E ^{\prime }] \mid p ^{n\mu }$. Observing
finally that $[\Delta ^{\prime }] \in {\rm Br}(Y ^{\prime }/E
^{\prime })$ (cf. \cite{P}, Sects. 9.4 and 13.1), one obtains that
ind$(\Delta ^{\prime }) \mid [Y ^{\prime }\colon E ^{\prime }] \mid
p ^{n\mu }$, so Lemma \ref{lemm4.1} is proved.
\end{proof}

\medskip
\begin{rema}
\label{rema4.2} Note that a field $E$ satisfies abrd$_{p}(E) <
\infty $, for some $p \in \mathbb P$, if and only if there exists $c
_{p}(E) \in \mathbb N$, such that each $A _{R} \in s(R)$ with exp$(A
_{R}) = p$ is Brauer equivalent to a tensor product of $c _{p}(E)$
algebras from $s(R)$ of degree $p$, where $R$ ranges over Fe$(E
_{p})$ and $E _{p}$ is the fixed field of a Sylow pro-$p$-subgroup
$G _{p}$ of $\mathcal{G}_{E}$. Since $E _{p}$ contains a primitive
$p$-th root of unity unless $p = {\rm char}(E)$, this can be deduced
from Lemma \ref{lemm4.1} and "quantative" versions of \cite{MS},
(16.1), and \cite{A1}, Ch. VII, Theorem~28 (see \cite{Ti}, page 506,
and \cite{Te}, respectively). When abrd$_{p}(E) < \infty $ and $p
\neq {\rm char}(E)$, $c _{p}(E)$ is in fact a cohomological
invariant of $G _{p}$ (cf. \cite{MS}, (11.5)). As noted in
\cite{Ka}, the Bloch-Kato Conjecture, proved in \cite{Voe}, implies
that if abrd$_{p}(E) < \infty $, then cd$_{p}(\mathcal{G}_{E}) <
\infty $ unless $E$ is formally real and $p = 2$ (see also \cite{L},
Ch. XI, Sect. 2, and \cite{Se}, Ch. I, 3.3).
\end{rema}

\medskip
Lemma \ref{lemm3.1} and the following two lemmas form the 
valuation-theoretic basis for the proof of the main results of this 
paper.

\medskip
\begin{lemm}
\label{lemm4.3} Let $(K, v)$ be a valued field with {\rm char}$(K) =
q > 0$, and let $K _{v}$ be a Henselization of $K$ in $K _{\rm sep}$
relative to $v$. Then:
\par
{\rm (a)} {\rm Brd}$_{q}(K) \le n$, provided that $[K\colon K ^{q}]
= q ^{n} < \infty $;
\par
{\rm (b)} {\rm Brd}$_{q}(K) \ge n$, if $v(K)/qv(K)$ has order $q
^{n}$ and $r _{q}(\widehat K) \ge n$; in this case, $(q ^{n}, q)$ is
an index-exponent pair over $K$.
\end{lemm}

\medskip
\begin{proof}
Lemma \ref{lemm4.3} (a) follows from (1.3), Lemma \ref{lemm4.1}, and
\cite{A1}, Ch. VII, Theorem~28, so it remains for us to prove Lemma
\ref{lemm4.3} (b). It is clear from the Artin-Schreier theorem that
$K$ possesses degree $q$ extensions $U _{1}, \dots , U _{n}$ in
$K(q)$, such that $U _{j} ^{\prime } = U _{j}K _{v}$ is inertial
over $K _{v}$ with $[U _{j} ^{\prime }\colon K _{v}] = q$, $j = 1,
\dots , n$, and $[U ^{\prime }\colon K _{v}] = q ^{n}$, where $U
^{\prime } = U _{1} ^{\prime } \dots U _{n} ^{\prime }$. As
$v(K)/qv(K)$ is of order $q ^{n}$, this enables one to deduce Lemma
\ref{lemm4.3} (b) from (3.6) (a).
\end{proof}

\medskip
\begin{lemm}
\label{lemm4.4} Let $(K, v)$ be a Henselian field with
$\mathcal{G}_{\widehat K}$ pronilpotent and {\rm
cd}$_{p}(\mathcal{G}_{\widehat K}) \le 1$, for some $p \in \mathbb
P$, $p \neq {\rm char}(\widehat K)$. Let also $\tau (p)$ be the
$\mathbb F _{p}$-dimension of $v(K)/pv(K)$, $\varepsilon _{p} \in
\widehat K _{\rm sep}$ a primitive $p$-th root of unity, and $m _{p}
= {\rm min}\{\tau (p), r _{p}(\widehat K)\}$. Then:
\par
{\rm (a)} {\rm Brd}$_{p}(K) = \infty $ if and only if $m _{p} =
\infty $ or $\tau (p) = \infty $ and $\varepsilon _{p} \in \widehat
K$; {\rm abrd}$_{p}(K) = \infty $, if and only if $\tau (p) = \infty
$;
\par
{\rm (b)} {\rm Brd}$_{p}(K) = {\rm abrd}_{p}(K) = [(m _{p} + \tau
(p))/2]$, in case $\varepsilon _{p} \in \widehat K$, $r
_{p}(\widehat K) \le 1$ and $\tau (p) < \infty $; {\rm Brd}$_{p}(K)
= m _{p}$, if $\varepsilon _{p} \notin \widehat K$ and $m _{p} <
\infty $.
\par
{\rm (c)} {\rm abrd}$_{p}(K) = \tau (p)$, if $r _{p}(\widehat K) \ge
2$ and $\tau (p) < \infty $.
\end{lemm}

\medskip
\begin{proof}
It is clear from (3.3) and (3.6) (a) that Brd$_{p}(K) = \infty $,
provided that $m _{p} = \infty $. Henceforth, we assume that $m _{p}
< \infty $. Suppose first that $\varepsilon _{p} \notin \widehat K$.
Since $p \neq {\rm char}(K)$ and cd$_{p}(\mathcal{G}_{\widehat K})
\le 1$ (whence, abrd$_{p}(\widehat K) = 0$), the Henselity of $v$
ensures that every $D \in d(K)$ with $[D] \in {\rm Br}(K) _{p}$ is a
nicely semi-ramified algebra over $K$, in the sense of \cite{JW}
(see Lemmas~5.14 and 6.2 therein). Hence, by \cite{JW}, Theorem~4.4,
$D$ is defined, for some $n \in \mathbb N$, by cyclic $p$-extensions
$U _{1}, \dots , U _{n}$ of $K$ in $K _{\rm ur}$, and by elements
$\pi _{1}, \dots , \pi _{n} \in K ^{\ast }$, in accordance with
(3.6) (a). This indicates that $U = U _{1} \dots U _{n}$ is a Galois
extension of $K$, $\mathcal{G}(U/K)$ has a system of at most $n$
generators, $n \le m _{p}$, and ind$(D)$ and exp$(D)$ are equal to
the order and to the period of $\mathcal{G}(U/K)$, respectively.
These observations prove that Brd$_{p}(K) = m _{p}$. Since
$[\widehat K(\varepsilon _{p})\colon \widehat K] \mid p - 1$, they
imply in conjunction with (1.1) (c), (3.2) and (3.3) that one may
assume, for the rest of the proof of Lemma \ref{lemm4.4}, that
$\varepsilon _{p} \in \widehat K$. Then (3.6) (b) yields Brd$_{p}(K)
= \infty $, if $\tau (p) = \infty $, so it remains for us to
consider the case of $\tau (p) < \infty $. As $p \neq {\rm
char}(\widehat K)$ and cd$_{p}(\mathcal{G}_{\widehat K}) \le 1$, it
is clear from (3.5) (b) and \cite{JW}, Lemmas~5.14 and 6.2, that
abrd$_{p}(K) = 0$, provided that $\tau (p) = 0$. This agrees with
the conclusions of the lemma, so we assume further that $\tau (p) >
0$. Our proof relies on the following observations:
\par
\medskip
(4.1) (a) For each $D \in d(K)$ with exp$(D) = p$, the group
$v(D)/v(K)$ has period $p$, $\widehat D$ is a field and $\widehat
D/\widehat K$ is a Galois extension, such that $\mathcal{G}(\widehat
D/\widehat K)$ is a homomorphic image of $v(D)/v(K)$; hence, ind$(D)
^{2} = [\widehat D\colon \widehat K]e(D/K) \mid p ^{m _{p}+\tau
(p)}$.
\par
(b) If $r _{p}(\widehat K) \ge 2$, then there exists a finite
extension $U$ of $K$ in $K _{\rm ur} \cap K(p)$, such that $r _{p}(U)
> \tau (p)$.
\par
\medskip\noindent
The inequality cd$_{p}(\mathcal{G}_{\widehat K}) \le 1$ ensures that
$\mathcal{G}(\widehat K(p)/\widehat K)$ is a free pro-$p$-group, so
(4.1) (b) can be deduced from (3.3) and Nielsen-Schreier's formula
for open subgroups of free pro-$p$-groups (cf. \cite{Se}, Ch. I,
4.2, and Ch. II, 2.1). Statement (4.1) (a) is contained in
\cite{JW}, (1.6) and Corollary~6.10). It follows from (4.1) (a) and
Lemma \ref{lemm4.1} that abrd$_{p}(K) \le \tau (p)$, whereas (3.6)
(a) and (4.1) (b) imply Brd$_{p}(U) \ge \tau (p)$. These results
prove Lemma \ref{lemm4.4} (a) and (c).
\par
We turn to the proof of Lemma \ref{lemm4.4} (b), so we assume that
$r _{p}(\widehat K) \le 1$. Then it follows from (3.2), (3.3),
\cite{Wh}, Theorem~2, and the conditions on $\mathcal{G}_{\widehat
K}$ that $r _{p}(\widehat K ^{\prime }) = r _{p}(\widehat K)$ and
$v(K ^{\prime })/pv(K ^{\prime }) \cong v(K)/pv(K)$, for every $K
^{\prime } \in {\rm Fe}(K)$. Hence, by (4.1) (a) and Lemma
\ref{lemm4.1}, Brd$_{p}(K ^{\prime }) \le [(m _{p} + \tau (p))/2]$,
proving that abrd$_{p}(K) \le [(m _{p} + \tau (p))/2]$. On the other
hand, it is clear from (3.6) (b) that Brd$_{p}(K) \ge [\tau (p)/2]$.
These observations prove Lemma \ref{lemm4.4} (b) in case $r
_{p}(\widehat K) = 0$. Suppose finally that $r _{p}(\widehat K) =
1$. Then $d(K)$ contains an algebra $B \otimes _{K} T$, defined in
accordance with (3.6) (c), for $n = 1$, $[U _{1}\colon K] = p$ and
$m = [(\tau (p) - 1)/2]$. In particular, ind$(B \otimes _{K} T) = p
^{m+1}$ and exp$(B \otimes _{K} T) = p$, which implies Brd$_{p}(K)
\ge 1 + m = [(1 + \tau (p))/2]$ and so completes our proof.
\end{proof}

\medskip
{\it We are now in a position to prove Theorem \ref{theo2.3}}. Let
$G$ be a pronilpotent group with cd$(G) = 1$, $G _{p}$ the Sylow
pro-$p$-subgroup of $G$ and $r _{p}$ the rank of $G _{p}$, for each
$p \in \mathbb P$. Suppose that $G _{2} \cong \mathbb Z _{2}$ and
put $\mathbb P ^{\prime } = \mathbb P \setminus \{2\}$. Then, it
follows from Burnside-Wielandt's theorem (cf. \cite{KM}, Ch. 6,
Theorem~17.1.4) that $G$ is isomorphic to the topological group
product $\prod _{p \in \mathbb P} G _{p}$. As cd$(G) = 1$, Lemma
\ref{lemm3.3} and Lemma \ref{lemm3.4}, applied to $P = \{2\}$, $G
_{2}$ and $\prod _{p \in \mathbb P'} G _{p}$, imply that $G \cong
\mathcal{G}_{K _{0}}$, for some characteristic zero field $K _{0}$
not containing a primitive $p$-th root of unity, for any $p \in
\mathbb P ^{\prime }$. This ensures that $r _{p}(K _{0}) = r _{p}$,
$p \in \mathbb P$. Note finally that $G$ can be chosen so that $r
_{p} = b _{p}$, $p \in \mathbb P ^{\prime }$, and by Lemma
\ref{lemm3.1}, there is a Henselian field $(K, v)$ with $\widehat K
\cong K _{0}$ and $n(p) = a _{p}$, $p \in \mathbb P ^{\prime }$.
Therefore, by Lemma \ref{lemm4.4}, we have Brd$_{p}(K) = b _{p}$ and
abrd$_{p}(K) = a _{p}$, for each $p \in \mathbb P ^{\prime }$.
Moreover, it follows from Lemmas \ref{lemm3.1} and \ref{lemm4.4}
that $(K, v)$ (specifically, the $\mathbb F _{2}$-dimension of
$v(K)/2v(K)$) can be chosen so that Brd$_{2}(K) = {\rm abrd}_{2}(K)
= a _{2}$. Theorem \ref{theo2.3} is proved.

\medskip
Our objective now is to prove Theorem \ref{theo2.2} in the case of
$q > 0$. The former part of this theorem is proved by applying our
next result to the field $K _{0} = \mathbb F _{q}$ and a system $c
_{p}$, $p \in \mathbb P$, with $c _{p} = \infty $, for all $p
\dagger q ^{2} - q$.

\medskip
\begin{lemm}
\label{lemm4.5} Let $K _{0}$ be a finite field with $q ^{m}$
elements, where $q = {\rm char}(K _{0})$. Put $P _{q,m} = \{p \in
\mathbb P\colon p \mid q(q ^{m} - 1)\}$, and fix a system $c _{p}
\in \mathbb N _{\infty }\colon \ p \in \mathbb P$. Then:
\par
{\rm (a)} There exists a Henselian field $(K, v)$ with {\rm
char}$(K) = q$, $\widehat K = K _{0}$, {\rm Brd}$_{q}(K) = c _{q}$
and {\rm abrd}$_{p}(K) = c _{p}$, for each $p \in \mathbb P$; this
ensures that {\rm Brd}$_{p}(K) \le 1$ when $p \in \mathbb P
\setminus P _{q,m}$, {\rm Brd}$_{p}(K) = c _{p}$, $p \in \mathbb P
_{q,m}$, and {\rm Brd}$_{p}(K) \neq 0$ in case $c _{p} \neq 0$;
\par
{\rm (b)} If $0 < c _{q} \neq \infty $, then $(K, v)$ can be chosen
so  that $[K\colon K ^{q}] = q ^{c _{q}}$;
\par
{\rm (c)} When $c _{q} = 0$, $(K, v)$ can be chosen so that $r
_{q}(K) = \infty $ and $K$ be perfect.
\end{lemm}

\medskip
\begin{proof}
Let $\bar n = n(p) \in \mathbb N _{\infty }\colon \ p \in \mathbb
P$, be a sequence, such that $n(p) = \infty $, provided $c(p) =
\infty $, and $2c _{p} - 1 \le n(p) \le 2c _{p}$ in case $c(p) <
\infty $ and $p \neq q$. Let also $(K, v)$ be a Henselian field with
char$(K) = q$ and $\widehat K = K _{0}$, attached to $\bar n$ as in
Lemma \ref{lemm3.1} (and subject to its additional restrictions in
case $c(q) < \infty $). Then it follows from Lemmas \ref{lemm4.3},
\ref{lemm4.4} and the equalities $r _{p}(K _{0}) = 1$, $p \in
\mathbb P$, that $(K, v)$ has the properties required by Lemma
\ref{lemm4.5}.
\end{proof}

\medskip
The extension $\Theta _{n}/\Lambda _{0}$ considered in Remark
\ref{rema3.2} has transcendency degree
\par\noindent
trd$(\Theta _{n}/\Lambda _{0}) = \infty $ (see \cite{BKu}, and
further references there). Hence, $\Lambda _{0}$ has a rational
extension $\Lambda _{\infty }$ in $\Theta _{n}$ with trd$(\Lambda
_{\infty }/\Lambda _{0}) = \infty $. This implies $[\Lambda \colon
\Lambda ^{q}] = [\Lambda _{\infty }\colon \Lambda _{\infty }^{q}] =
\infty $, where $\Lambda $ is the separable closure of $\Lambda
_{\infty }$ in $\Theta _{n}$. Therefore, the latter assertion of
Theorem \ref{theo2.2} can be deduced from Lemma \ref{lemm4.5} and
the following lemma.

\medskip
\begin{lemm}
\label{lemm4.6} Let $K _{0}$ be a finite field, and in the setting
of Remark \ref{rema3.2}, put $\Theta = \Theta _{n}$, and suppose
that $\Lambda \in I(\Theta /\Lambda _{0})$ is separably closed in
$\Theta $. Then:
\par
{\rm (a)} The valuations $w _{\Lambda }$, $\kappa _{\Lambda }$ and
$\theta _{\Lambda }$ of $\Lambda $ are Henselian;
\par
{\rm (b)} For each finite separable extension $R$ of $\Lambda $ in
$K _{\rm sep}$, $R\Theta $ is a completion of $R$ relative to the
topology induced by $w _{R}$, and $w _{R\Theta }$ is the continuous
prolongation of $w _{R}$ on $R\Theta $; in addition, $D _{R} \otimes
_{R} R\Theta \in d(R\Theta )$, for every $D _{R} \in d(R)$;
\par
{\rm (c)} The field $\Phi = \Lambda (\Gamma )$ satisfies the
equalities {\rm Brd}$_{p}(\Phi ) = {\rm Brd}_{p}(K)$ and
abrd$_{p}(\Phi ) = {\rm abrd}_{p}(K)$, $p \in \mathbb P$, {\rm
Brd}$_{q}(\Phi ) = {\rm abrd}_{q}(\Phi ) = n$, and $[\Phi \colon
\Phi ^{q}] = [\Lambda \colon \Lambda ^{q}]$.
\end{lemm}

\medskip
\begin{proof}
Lemma \ref{lemm4.6} (a) follows from \cite{Efr2}, Theorem~15.3.5,
and the Henselity of the valuations $w$, $\kappa $ and $\theta $ of
$\Theta $. The former claim of Lemma \ref{lemm4.6} (b) is obvious,
and it enables one to deduce the latter part of Lemma \ref{lemm4.6}
(b) from \cite{CohnPM}, Theorem~2. As $v = w _{K}$, $w _{\Phi }(\Phi
) = w _{K}(K)$ and $K _{0}$ is the residue field of $(K, v)$ and
$(\Phi , w _{\Phi })$, Lemma \ref{lemm4.4} implies Brd$_{p}(\Phi ) =
{\rm Brd}_{p}(K)$ and abrd$_{p}(\Phi ) = {\rm abrd}_{p}(K)$, for
each $p \neq q$. Observing that $[\Theta \colon \Theta ^{q}] = q
^{n}$, one obtains from Lemma \ref{lemm4.6} (b) and \cite{A1}, Ch.
VII, Theorem~28, that Brd$_{q}(R) \le {\rm Brd}_{q}(R\Theta ) \le
n$, for every finite separable extension $R$ of $\Lambda $ in $K
_{\rm sep}$. This proves that abrd$_{q}(\Lambda ) \le {\rm
abrd}_{q}(\Theta ) = n$, which leads to the conclusion that
Brd$_{q}(\Phi ) \le {\rm abrd}_{q}(\Phi ) \le {\rm abrd}_{q}(\Lambda
)$ (see also \cite{Ch2}, (1.3)). On the other hand, by Remark
\ref{rema3.2}, $\kappa _{\Phi }(\Phi ) = \kappa (K)$ and the residue
field of $(\Phi , \kappa _{\Phi })$ is isomorphic to $\Theta _{0}$.
Since, by the proof of Lemma \ref{lemm3.1}, $r _{q}(\Theta _{0}) =
\infty $ and $\kappa (K)/q\kappa (K)$ is of order $q ^{n}$, this
allows us to obtain from Lemma \ref{lemm4.3} that Brd$_{q}(\Phi )
\ge n$. Note finally that $\Phi /\Lambda $ is a separable extension,
so $[\Phi \colon \Phi ^{q}] = [\Lambda \colon \Lambda ^{q}]$, which
completes our proof.
\end{proof}

\medskip
\begin{rema}
\label{rema4.7} The proof of Theorem \ref{theo2.2} is technically
simpler in characteristic $2$. Lemma \ref{lemm4.4} shows that if $K
_{0} = \mathbb F _{2}$ and $\Theta _{0}$ is a perfect closure of the
extension $K _{\infty }$ of $K _{0}$ defined in the proof of Lemma
\ref{lemm3.1}, then abrd$_{2}(\Theta _{0}) = 0$, Brd$_{p}(\Theta
_{0}) = 1$ and abrd$_{p}(\Theta _{0}) = \infty $, for all $p > 2$.
When $n \in \mathbb N$, $\Theta _{n}$ and $\Lambda _{0}$ are defined
as in Remark \ref{rema3.2}, $\Lambda _{\infty }$ is a rational
extension of $\Lambda _{0}$ in $\Theta _{n}$ with trd$(\Lambda
_{\infty }/\Lambda _{0}) = \infty $, and $\Lambda $ is the separable
closure of $\Lambda _{\infty }$ in $\Theta _{n}$, we have $[\Lambda \colon \Lambda ^{2}] = \infty $, Brd$_{2}(\Lambda ) = {\rm
abrd}_{2}(\Lambda ) = n$, and for each $p > 2$, Brd$_{p}(\Lambda ) =
1$ and abrd$_{p}(\Lambda ) = \infty $. Note also, omitting the
details, that the field $\Theta _{0}$ enables one to find an
alternative proof of Theorem \ref{theo2.2} in zero characteristic
(see \cite{Ch4}, Example~6.2).
\end{rema}

\medskip
When $c _{p} \in \mathbb N$, $p \in \mathbb P$, is an unbounded
sequence, the fields $E$ singled out by Lemma \ref{lemm4.5} have the
properties required by (2.4) (a). As to (2.4) (b), it is implied by
Lemma \ref{lemm3.1} and our next result.

\medskip
\begin{coro}
\label{coro4.8} In the setting of Lemma \ref{lemm4.4}, let $\widehat
K$ be a quasifinite field with {\rm char}$(\widehat K) = 0$ and
$\varepsilon _{p} \notin \widehat K$, for any $p \in \mathbb P
\setminus \{2\}$, and let $U _{n}$ be the degree $n$ extension of
$K$ in $K _{\rm ur}$, for a fixed integer $n \ge 2$. Suppose that $P
_{n} = \{p _{n} \in \mathbb P\colon \ n \mid p _{n} - 1\}$,
$[\widehat K(\varepsilon _{p _{n}})\colon \widehat K] = n$, for all
$p _{n} \in \mathbb P _{n}$, and the sequence $\tau (p)\colon \ p
\in \mathbb P$, satisfies the condition $\tau (p) = \infty $ if and
only if $p \in \mathbb P _{n}$. Then a field $L \in {\rm Fe}(K)$
satisfies ${\rm Brd}_{p}(L) < \infty $, $p \in \mathbb P$, if and
only if $U _{n} \notin I(L/K)$. When $U _{n} \notin I(L/K)$ and the
system $\tau (p)$, $p \in \mathbb P \setminus P _{n}$, is bounded,
{\rm Brd}$(L) < \infty $.
\end{coro}

\medskip
\begin{proof}
Lemma \ref{lemm4.4} and our assumptions show that if $p \notin P
_{n}$, then Brd$_{p}(L) \le {\rm abrd}_{p}(K) < \infty $. When $p
\in P _{n}$ and $L \in {\rm Fe}(K)$, they prove that Brd$_{p}(L) =
\infty $ if and only if $\varepsilon _{p} \in \widehat L$, and this
occurs if and only if $U _{n} \subseteq L$. The concluding assertion
of Corollary \ref{coro4.8} follows from Lemma \ref{lemm4.4}.
\end{proof}

\medskip
Lemmas \ref{lemm3.3} and \ref{lemm3.4} indicate that there exists a
quasifinite field $E$ of zero characteristic, such that
$[E(\varepsilon )\colon E] = 2 ^{y(p)}$, $p \in \mathbb P$, where
$\varepsilon _{p}$ is a primitive $p$-th root of unity in $E _{\rm
sep}$ and $y(p)$ is defined as in Lemma \ref{lemm3.3}, for each $p$.
At the same time, Lemma \ref{lemm3.1} and Corollary \ref{coro4.8}
imply the existence of Henselian fields $(E _{n}, v _{n})$ with
$\widehat E _{n} = E$, which possess the properties required by
(2.4) (b), for $n = 2 ^{t}$, $t \in \mathbb N$. Using \cite{Ch3},
Lemma~3.2, instead of Lemma \ref{lemm3.3}, and arguing in the same
way, one proves (2.4) (b) in general.

\vskip0.5truecm\noindent \emph{Acknowledgements.} A considerable
part of this research was carried out during my visit to Tokai 
University, Hiratsuka, Japan, in 2012. I would like to thank my 
host-professor Junzo Watanabe, the colleagues at the Department of 
Mathematics, and Mrs. Yoko Kinoshita and her team for their genuine 
hospitality. I also acknowledge that my research was partially 
supported by a Project No. RD-08-241/12.03. 2015 of Shumen 
University, Bulgaria.

\par


\medskip
\begin{thebibliography}{aa}

\bibitem{A1} A.A. Albert, \emph{Structure of Algebras}, Amer. Math.
Soc. Colloq. Publ., vol. XXIV, 1939.

\bibitem{ABGV} A. Auel, E. Brussel, S. Garibaldi, U. Vishne,
\emph{Open problems on central simple algebras}, Transform. Groups
{\bf 16} (2011), 219-264.

\bibitem{BKu} A. Blaszczok, F.-V. Kuhlmann, \emph{Algebraic
independence of elements in immediate extensions of valued fields},
Preprint, arXiv:1304.1381v1 [math.AC].

\bibitem{Ch1} I.D. Chipchakov, \emph{The normality of locally finite
associative division algebras over classical fields}, Vestn. Mosk.
Univ., Ser. I (1988), No. 2, 15-17 (Russian: English transl. in:
Mosc. Univ. Math. Bull. {\bf 43} (1988), 2, 18-21).

\bibitem{Ch2} I.D. Chipchakov, \emph{On the residue fields of
Henselian valued stable fields}, J. Algebra {\bf 319} (2008), 16-49.

\bibitem{Ch3} I.D. Chipchakov, \emph{On Brauer $p$-dimensions and
absolute Brauer $p$-dimensions of Henselian fields}, Preprint,
arXiv:1207.7120v4 [math.RA].

\bibitem{Ch4} I.D. Chipchakov, \emph{On Brauer $p$-dimensions and
index-exponent relations over finitely-generated field extensions},
Preprint, arXiv:1501.05977v1 [math.RA].

\bibitem{CohnPM} P.M. Cohn, \emph{On extending valuations in
division algebras}, Stud. Sci. Math. Hung. {\bf 16} (1981), 65-70.

\bibitem{Dr1} P.K. Draxl, \emph{Ostrowski's theorem for Henselian
valued skew fields}, J. Reine Angew. Math. {\bf 354} (1984),
213-218.

\bibitem{Efr2} I. Efrat, \emph{Valuations, Orderings, and Milnor
$K$-Theory}, Math. Surveys and Monographs, 124, Providence, RI:
Amer. Math. Soc., XIII, 2006.

\bibitem{FV} I.B. Fesenko, S.V. Vostokov, \emph{Local Fields and
Their Extensions}, 2nd ed., Transl. Math. Monographs, 121, Amer.
Math. Soc., Providence, RI, 2002.

\bibitem{HHKr} D. Harbater, J. Hartmann, D. Krashen,
\emph{Applications of patching to quadratic forms and central simple
algebras}, Invent. Math. {\bf 178} (2009), 231-263.

\bibitem{JW} B. Jacob, A. Wadsworth, \emph{Division algebras
over Henselian fields}, J. Algebra {\bf 128} (1990), 126-179.

\bibitem{Jong} A.J. de Jong, \emph{The period-index problem for the
Brauer group of an algebraic surface}, Duke Math. J. {\bf 123}
(2004), 71-94.

\bibitem{Ka} B. Kahn, \emph{Comparison of some field invariants},
J. Algebra {\bf 232} (2000), 485-492.

\bibitem{KM} M.I. Kargapolov, Yu.I. Merzlyakov,
\emph{Fundamentals of Group Theory}, 3rd Ed., Nauka, Moscow, 1982.

\bibitem{L} S. Lang, \emph{Algebra}, Addison-Wesley Publ. Comp.,
Mass., 1965.

\bibitem{Lieb} M. Lieblich, \emph{Twisted sheaves and the
period-index problem}, Compos. Math. {\bf 144} (2008), 1-31.

\bibitem{Matz} E. Matzri, \emph{Symbol length in the Brauer group of
a field}, Preprint, arXiv:1402.0332v1 [math.RA].

\bibitem{MS} A.S. Merkur'ev, A.A. Suslin, \emph{$K$-cohomology of
Severi-Brauer varieties and norm residue homomomorphisms}, Izv.
Akad. Nauk SSSR {\bf 46} (1982), 1011-1046 (Russian: English transl.
in: Math. USSR Izv. {\bf 21} (1983), 307-340).

\bibitem{Mo} P. Morandi, \emph{The Henselization of a valued
division algebra}, J. Algebra {\bf 122} (1989), 232-243.

\bibitem{PaSur} R. Parimala, V. Suresh, \emph{Period-index and
$u$-invariant questions for function fields over complete discretely
valued fields}, Preprint, arXiv:1304.2214v1 [math.RA].

\bibitem{P} R. Pierce, \emph{Associative Algebras}, Graduate Texts
in Math., vol. 88, Springer-Verlag, XII, New York-Heidelberg-Berlin,
1982.

\bibitem{Re} M. Reiner, \emph{Maximal Orders}, London Math. Soc.
Monographs, vol. 5, London-New York-San Francisco: Academic Press, a
subsidiary of Harcourt Brace Jovanovich, Publishers, 1975.

\bibitem{Sch} O.F.G. Schilling, \emph{The Theory of Valuations},
Mathematical Surveys, No. 4, Amer. Math. Soc., New York, N.Y., 1950.

\bibitem{Se} J.-P. Serre, \emph{Galois Cohomology}, Transl. from the
French original by Patrick Ion, Springer, Berlin, 1997.

\bibitem{Te} O. Teichm\" uller, \emph{Zerfallende zyklische
$p$-algebren}, J. Reine Angew. Math. {\bf 176} (1937), 157-160.

\bibitem{Ti} J.-P. Tignol, \emph{On the length of decompositions of
central simple algebras in tensor products of symbols}, in: Methods
in Ring Theory, Proc. NATO Adv. Study Inst., Antwerp/Belg. 1983,
NATO ASI Ser., Ser. C {\bf 129} (1984), 505-516.

\bibitem{Voe} V. Voevodsky, \emph{On motivic cohomology with
$\mathbb Z/l$-coefficients}, Ann. Math. {\bf 174} (2011), 401-438.

\bibitem{W} S. Warner, \emph{Topological Fields}, North-Holland
Math. Studies, 157; Notas de Mat\'{e}matica, 126. North-Holland
Publishing Co., Amsterdam, 1989.

\bibitem{Wat} W.C. Waterhouse, \emph{Profinite groups are Galois
groups}, Proc. Amer. Math. Soc. {\bf 42} (1974), 639-640.

\bibitem{Wh} G. Whaples, \emph{Algebraic extensions of arbitrary
fields}, Duke Math. J. {\bf 24} (1957), 201-204.

\end{thebibliography}
\end{document}